\documentstyle[12pt]{article}
\newcommand{\f}{\mbox{$\cal F$~}}

\newcommand{\arrow}{\rightarrow}
\newcommand{\ar}{\rightarrow}
\newcommand{\var}{\mbox{$\vartheta$}}

\newcommand{\atla}{\mbox{} \newline \mbox{} \newline}
\newcommand{\qed}{\hfill \mbox{$\Box$} \\ \mbox{} \\ }

\newcommand{\C}{{\bf C}}

\newcommand{\ben}{\begin{eqnarray*}}
\newcommand{\een}{\end{eqnarray*}}
\newcommand{\be}{\begin{eqnarray}}
\newcommand{\ee}{\end{eqnarray}}
\newcommand{\beg}{\begin{equation}}
\newcommand{\eg}{\end{equation}}

\newcommand{\del}{\partial}
\begin{document}
\mbox{}
\\
\begin{center}
{\bf GENERIC SINGULARITIES OF HOLOMORPHIC FOLIATIONS } \\
Sinan Sert\"{o}z \\
Bilkent University
\end{center}
{}~~~ \\
{}~~~~
\par
Singular behaviour of holomorphic dynamics has been of great interest
to researchers of related areas.
Baum has analyzed the behaviour of leaves
of a holomorphic foliation
around a singularity where the singular set has one less dimension
than the leaf dimension,  \cite{BA}.
He then used this with
Bott to describe their by now well known residue classes, under the same
requirement about the dimension of the singular set,
\cite{BB}, see also
\cite{VI}.
Later Cenkl showed
that Baum-Bott residue classes can be described also for foliations
where the singular set has smaller dimension, \cite{CE}.
This brought up the question
about the structure of such singularities.
In this article we study the structure of the singularity of a
holomorphic foliation
for which the dimension of the singular set is not necessarily one less
than the leaf rank. We show that
the structure of the leaves around a
point on the singular set is determined by
the vector space rank of the sheaf defining the foliation at this
point. The generic singularity is then the
one for which this rank is zero.
We also show that in general the singularity of a holomorphic foliation
is locally the pullback of a generic singularity.
\par
In this article we also give examples of holomorphic foliations of
rank $k$, on an $n$
dimensional manifold, with
a singular set of dimension $r$, where $r$ is any natural number less
than $k$, subject only to the condition that $r\geq n-k-1$ which
is exactly a condition stated by Malgrange in \cite{MA}.
Such examples are missing in the literature.
The examples given here fill in
this gap and also illustrate the claims of our main theorem.
\par
We have utilized Baum's techniques of \cite{BA} where, with the above
notation, he has studied the case $r=k-1$ which is also the
case studied by Vishik in \cite{VI}.
\par
Many authors such as Verjovsky, Alexander and Thom has analyzed foliations of
rank $k$ with isolated singularities, which justifies our results,
\cite{AV},\cite{GSV}, \cite{TH2}.
See also \cite{XA} where foliations with curves having isolated singularities
are examined.
\par
I thank Professor A. Verjovsky  for his numerous comments
on this problem during my visit to Trieste. I also thank Professor
J.B. Carrell who read an earlier version of this article. His
suggestions improved the article considerably. \\
{\bf The set up and the problem } \par
Throughout the article $M$ will
denote a complex manifold of dimension $n$.
$TM$ will denote its holomorphic tangent bundle and $\tau$ the sheaf of
holomorphic sections of $TM$. We will use $\xi$ to denote an integrable
full coherent subsheaf of $\tau$ of fixed rank $k$.
(For the definition of ``full"
see  \cite{BB}.) $Z$ will denote a fixed connected component
of the singular set $S$ of $\xi$. We will
assume that $Z$ is irreducible of fixed dimension $r$.
{\var}  will
denote the structure sheaf of $M$ and $\mu$ will be the maximal ideal of
{\var}.
\par
For each point $x\in M$, $T_x(\xi )$ will denote the vector
subspace of $T_xM$ defined as $\xi_x/\xi_x\otimes \mu_x$. For an
equivalent definition of $T_x(\xi )$ using vector fields see either
\cite{BB} or \cite{BA}.
\par
There is a natural filtration of $Z$
depending on the behaviour of $T_x(\xi )$
on $Z$;
\[ Z=Z^1 \supset Z^2 \supset\cdots\supset Z^k \]
where
\[ Z^i=\{ x\in Z \; | \; \dim T_x(\xi ) \leq k-i\; \} . \]
Each $Z^i$ is a closed subvariety.
\par
To formulate the results it will be convenient to introduce the
following formalism: Let $Fol(n,k,r,s)$ denote the set of isomorphism
classes of singular holomorphic foliations around the origin in $\C^n$
defined as follows:
An element  in $Fol(n,k,r,s)$ is represented by an equivalence class
of pairs $<{\f},U>$ where \\
1) $U$ is an open neighbourhood of the origin in $\C^n$ and {\f} is an
integrable full coherent subsheaf of the tangent
sheaf of $\C^n$ on $U$. \\
2) rank{\f}=$k$. \\
3) The singular set $Z$ of {\f} is connected, smooth,
contains
the origin and is of dimension $r$. \\
4) The dimension of $T_x(\xi )$ is $s$ for all $x$ in $Z$.
i.e. we assume that
$Z=Z^{k-s}-Z^{k-s+1}$. This corresponds to locally choosing $U$ small
enough to exclude $Z^{k-s+1}$. \\
(For alternate definitions see \cite{AV}, \cite{SU2} and \cite{CA}). \\
Two elements $<\f_1,U_1>$ and $<\f_2,U_2>$ represent the same element
if there is an invertible
holomorphic map $f$ of $W\subset U_1\cap U_2$ onto $W$
with $f^{\ast}\f_1=\f_2$.
\par
We are interested in the cases when $n>k>r\geq s \geq 0$.
\par
With any element $\xi =<{\f},U>$ of $Fol(n,k,r,s)$ we use the notation
$\xi\times\C^m$
to denote that element of
$Fol(n+m,k+m,r+m,s+m)$
which represents the foliation on $\C^n\times\C^m$ given by the
geometric pull back of the leaves of $\xi$ on $\C^n$.
The notation $Fol(n,k,r,s)\times\C^m$ then denotes the set of all
$\xi\times\C^m$ where $\xi\in Fol(n,k,r,s)$.
\par
We will need the following map
\[ i_m \; : \; Fol(n,k,r,0)\arrow Fol(n+m,k,r+m,0) \]
where if $\alpha$ is in $Fol(n,k,r,0)$ and singularly foliates $\C^n$
with leaves
$\{ L\}$ then the leaves of $i_m(\alpha )$ on $\C^n\times\C^m$
are all of the form $L\times \{ t\}$ for $t\in\C^m$, i.e. each slice
$\C^n\times \{ t\}$ is foliated as prescribed by $\alpha$.
We denote this phenomena by
\[ i_m(\alpha )=\bigcup_{t\in\C^m}\alpha\times\{ t\}. \]
Now we can state our main result: \atla
{\bf Main Theorem: } $Fol(n,k,r,s)=Fol(n-s,k-s,r-s,0)\times\C^s$. \\
{\sf i.e. for every $\xi$ in $Fol(n,k,r,s)$ there is an $\eta_{\xi}$ in
$Fol(n-s,k-s,r-s,0)$ such that $\xi =\eta_{\xi}\times\C^s$. }
\atla
Furthermore we will call $\xi$ in $Fol(n,k,r,s)$ ``split" if
the corresponding $\eta_{\xi}$,
whose existence is established by the above
theorem , is of the form $i_{r-s}(\alpha )$ for some $\alpha$ in
$Fol(n-r,k-s,0,0)$.
A necessary condition for $\xi$ to be split is $n-r>k-s$.
A split foliation locally looks like a collection of foliations with
isolated singularities; see the examples' section.
\atla
{\bf The proof of the main theorem } \par
Let $p$ be a smooth point of $Z$ and choose a local coordinate system
$(z,U)$ centered at $p$ such that
\[ Z\cap U=\{ q\in U| z_{r+1}(q)=\cdots z_{n}(q)=0\;\}. \]
Baum's tangency lemma, \cite{BA},
forces $T_x(\xi )$ to stay in $T_x(Z)$ for
all $x\in Z$. Since the rank of $T_x(\xi )$ is fixed throughout $Z$ and
$\xi$ is coherent,
we can choose vector fields
$V_1,...,V_s$ on $U$
such that \\
1) $V_1(x),...,V_s(x)$ are in $T_x(\xi )$ for all $x\in U$. \\
2) $V_i(x)=\frac{\del}{\del z_i}|_x$ for $x\in Z$, $i=1,2,...,s.$ \\
3) $V_1(x),...,V_s(x), \frac{\del}{\del z_{s+1}}|_x,...,
\frac{\del}{\del z_n}|_x$ are linearly independent for all $x\in U$. \\
This suggests that we write $Z$ as $AZ\times BZ$ where
\ben
 AZ=\{ x\in Z|\; z_{s+1}(x)=\cdots =z_n(x)=0 \;\} \\
BZ=\{ x\in Z|\; z_1(x)=\cdots =z_s(x)=0
\; {\rm and }\; z_{r+1}(x)=\cdots
=z_n(x)=0 \; \} .
\een
Note that $\dim AZ=s$ and $\dim BZ=r-s$.
We also define a disc $D$ in $U$ as
\[ D=\{ x\in U|\; z_1(x)=\cdots =z_s(x)=0 \;\} . \]
Then the dimension of $D$ is $n-s$ and $U$
can be realized as $AZ\times D$
in the usual sense.
\par
The vector fields $V_1,...,V_s$
generate a locally free subsheaf $\zeta$ of
$\xi$ with rank $s$.
At each point of $U$ the complement of $\zeta$ in $\xi$
is a $k-s$ dimensional vector space, except on $Z$ where the complement
is zero. Thus $\xi$ can be realized as $\zeta\oplus\eta$, where
$\eta$ is a coherent subsheaf of $\xi$ with rank $k-s$ on $U-Z$ and
$0$ on $Z$. Note that since $\xi$ is coherent it is generated by
$V_1,...,V_s$ and some other sections $v_1,...,v_d$ around
$p$, and $\eta$ is then the sheaf generated by $v_1,...,v_d$.
\par
For all $x$ in $D$, the vector space corresponding to $\eta_x$,
(the vector space generated by germs of $\eta_x$ evaluated at $x$),
lies in the
space generated by
$\frac{\del}{\del z_{s+1}}|_x,...,\frac{\del}{\del z_n}
|_x$ because of the choice of $V_i$'s. Hence $\eta |D$ lies in the
tangent sheaf of $D$, and consequently so does $[\eta |D,\eta |D]$
where $[\cdot ,\cdot ]$ denotes the lie bracket.
Also $[\eta ,\eta ]\subset \xi$, and since $\xi |D$ and the
tangent sheaf of $D$ have only
$\eta |D$ in common we conclude that $\eta |D$
is integrally closed.
\par
Hence $\eta_{\xi}=\eta |D$ foliates $D$ with
leaves of rank $k-s$ and with a
singularity along $BZ$. Since $BZ$ is in $Z$ and the rank of $\eta$ is
zero on $Z$, the rank of $\eta_{\xi}$ is zero along its singularity.
\par
Thus $\eta_{\xi}$ is in $Fol(n-s,k-s,r-s,0)$ around $p$.
\par
Now we use Baum's trick to
carry the foliation of $D$ over $U$ along $AZ$,
\cite{BA}.  For this we define a vector field $V_t$
for each point $t=(t_1,...,t_s,0,...,0)$ of $AZ$;
\[ V_t(x)=t_1V_1(x)+\cdots +t_sV_s(x), \;\;\; x\in U. \]
We will move $D$ with the flow generated by $V_t$;
\[ D_t:=exp(V_t(D)). \]
We assume that $U$ is shrunk so that all the $D_t$'s will
lie in $U$,
and will be disjoint  for different $t$. Since $V_t(p)$ lies in
the tangent space of $AZ$,
$exp(V_t(p))$ will also lie in $AZ$. We expect
$U$ to be small so that
$exp(V_t(D))$ will not intersect $AZ$ in any point
other than $exp(V_t(p))$.
\par
We will now find the intersection point of $exp(V_t(D))$ with $AZ$.
Recall that $exp_{V_t(p)}:\C \ar U$ is the unique map
$y\mapsto (e_1(y),...,e_n(y))$ which sends $0$ to $p$ and whose
differential sends $\frac{d}{dy}$ of $T_0\C$ to $V_t(p)$ of $T_pU$.
In particular
\ben
(d\; exp_{V_t(p)})(\frac{d}{dy})(z_i)  & = &
       V_t(p)(z_i) \\
 & = & (t_1 \frac{\partial}{\partial z_1}|_p+\cdots
           +t_{k-r}\frac{\partial}{\partial z_s}|_p)(z_i) \\
 & = & \left\{ \begin{array}{ll}
      z_i & {\rm if }\;\; 0<i\leq s, \\
      0   & {\rm if }\;\; s<i\leq n.
      \end{array} \right.
\een
On the other hand
\ben
(d\; exp_{V_t(p)})(\frac{d}{dy})(z_i) & = &
                  \frac{d}{dy}(z_i\circ (exp_{V_t(p)}(y))) \\
 & = & \frac{d}{dy}(z_i \circ (e_1(y),...,e_n(y))) \\
 & = & \frac{d}{dy}e_i(y)
\een
Hence, noting that $0$ goes to $p$, we have
\[ exp_{V_t(p)}(y)=(yt_1,...,yt_s,0,...,0). \]
And finally,
\ben
exp(V_t(p)) & = & exp_{V_t(p)}(1) \\
            & = & (t_1,...,t_s,0,...,0) \\
            & = & t,
\een
thus showing that
\[ D_t\cap AZ=\{ t\}. \]
Since $V_t(x)$ is in $T_x(\xi )$ for all $x$ in $U$,  the flow passing
through $x$, the integral curve of $V_t$, remains in the same leaf of
$\xi$ as $x$. \\
For every $x\in U$ we can find a unique $t=t(x)$ in $AZ$ and  a unique
$q=q(x)$ in $D$, choosing $U$ smaller if necessary, such that
\[ x = exp(V_t(q)). \]
Since $t=(t_1(x),...,t_s(x),0,...,0)$ and
$q=(0,...,0,q_{s+1}(x),...,q_n(x))$ we can
choose to denote $x$ with the coordinates
$(t_1(x),...,t_s(x),q_{s+1}(x),...,q_n(x))$.
This makes $U$ holomorphically equivalent to $AZ\times D$ around $p$
with respect to these new coordinates,
which we agree to denote by $Z\odot D$. \\
Let $pr_2:U\ar D$ be the map which sees $U$ as $AZ\odot D$ and projects
to the second component;
\[ x\mapsto (q_{s+1}(x),...,q_n(x)). \]
If $D$ is foliated by $\eta_{\xi}$ with leaves $\{ L\}$, as mentioned
before, then any leaf  $L$ is $k-s$ dimensional and is carried
along the flow lines of $V_t$ with the parameter $t$ lying in an $s$
dimensional space $AZ$. Hence the images of $L$ along these flow lines,
all remaining in the same leaf of $\xi$ as $L$, fill out a $k$
dimensional space, which must in turn be a whole leaf of $\xi$ in $U$.
Thus $U=AZ\odot D$ is foliated by leaves of the form $L\times AZ\simeq
L\times\C^s$ with a singularity along $BZ$, which means finally that
$\xi$ is of the form $\eta_{\xi}\times\C^s$ where $\eta_{\xi}$ is in
$Fol(n-s,k-s,r-s,0)$.
\qed
\par
One interesting case is when the foliation $\eta_{\xi}$ can further be
decomposed; for this let $K$ be the disc
\[ K=\{ x\in U|\; z_1(x)=\cdots =z_r(x)=0\;\}, \;\;\; \dim K=n-r. \]
$K$ is a disc in $D$ perpendicular to $BZ$ at $p$.
If $K$ intersects the leaves of $\xi$ transversally, then the leaves of
$\xi$ singularly foliate
$K$ with $k-s$ dimensional leaves and with an isolated singularity.
This foliation can be moved along $AZ$ via $V_t$ as before to
foliate $K\times AZ$ with a rank $k$ foliation and with a
singularity along $AZ$.
We then have to repeat this success at all the other points of $BZ$,
besides $p$, (see example 3).
\par
\mbox{} \\
{\bf Examples } \\
{\bf 1) } The lines through $0$ foliate $\C^n\backslash 0$
thus showing
that $Fol(n,1,0,0)$ is not empty. Next consider the following foliation
by divisors: for any $\lambda\in \C$ let
\[ L_{\lambda}=\{ (X_1,...,X_n)\in\C^n|\;
X_1^2+\cdots +X_n^2=\lambda \;\} .\]
Each $L_{\lambda}$ is smooth except when $\lambda$ is zero, in
which case $L_{\lambda}$ has an isolated singularity at the origin.
The collection $\{ L_{\lambda}\}$ foliates $\C^n$ with an isolated
singularity at the origin. Thus $Fol(n,n-1,0,0)\neq\emptyset$. \\
{\bf 2) } For each $\lambda =(\lambda_1,\lambda_2)\in\C^2$ define
\[ L_{\lambda}=\{ (X_1,...,X_n)\in\C^n|\; X_1^2+\cdots +X_{n-1}^2
=\lambda_1, \;\; X_n=\lambda_2 \;\} .\]
Each $L_{\lambda}$ is smooth except when $\lambda_1=0$, in which case
$L_{(0,\lambda_2)}$ has an isolated singularity at
$(0,...,0,\lambda_2)$.
The set $\{ L_{\lambda}\}$ foliates $\C^n$ with a singularity along
the last coordinate. Thus $Fol(n,n-2,1,0)$ is not empty.
Since $s=0$ we can check if the foliation is split.
In this example $AZ$ is just the origin because
each $L_{(\lambda_1,\lambda_2)}$ lies
in the hyperplane $X_n=\lambda_2$, and the tangent vectors of the leaf
do not have an $X_n$ component.
The leaves hence intersect the hyperplane
transversally. This means that each such hyperplane is foliated as in
example~1.  Let the foliation of the hyperplane be $\eta$, which belongs
to $Fol(n-1,n-2,0,0)$. Then the foliation of this example is of the
form $i_1(\alpha )\times\C^1$ and is therefore split. \\
{\bf 3) } We now give a less obvious example of a split foliation. We will
foliate $\C^n$ with rank $k$ leaves and with a singularity along
an $r$ dimensional subvariety. For this type of construction to work
we must necessarily have $r\geq n-k-1$, which turns out to be a
bound given by Malgrange in \cite{MA}. \\
For each $\lambda =(\lambda_1,...,\lambda_{n-k})\in\C^{n-k}$
we define the leaf $L_{\lambda}$ as
\ben
L_{\lambda}=\{ (X_1,...,X_n)\in\C^n | & &
     X_1^2+\cdots +X_{n-r}^2=\lambda_1, \\
 & & X_{n-r+1}=\lambda_2, \\
 & & \vdots \\
 & & X_{2n-k-r-1}=\lambda_{n-k}. \;\;\}
\een
Dimension of each leaf is $k$. When $\lambda_1\neq 0$ then $L_{\lambda}$
is smooth and when $\lambda_1=0$ then $L_{\lambda}$ is singular along
a $k+r+1-n$ dimensional subvariety.
The set of singular leaves is $n-k-1$
dimensional. Hence the dimension of the singularity is $r$.
Let us call this foliation $\xi$. By examining the Jacobian of the
system giving the foliation we can find that $s=k+r+1-n$. Hence
\[ \xi\in Fol(n,k,r,k+r+1-n). \]
We note that
\ben
Z=\{ (0,...,0,X_{n-r+1},...,X_n)\in\C^n\} ,\;\;\dim Z=r. \\
AZ=\{ (0,...,0,X_{2n-k-r},...,X_n)\in\C^n\} ,\;\;\dim AZ=k+r+1-n. \\
BZ=\{ (0,...,0,X_{n-r+1},...,X_{2n-k-r-1},0,...,0)\in\C^n\} ,\;\;\dim BZ
                                                             =n-k-1. \\
D=\{ (X_1,...,X_{2n-k-r-1},0,...,0)\in\C^n\} ,\;\;\dim D=2n-k-r-1.
\een
$D$ is foliated by the same set of $n-k$ equations that foliate $\C^n$.
Call this foliation $\eta_{\xi}$. The singular set of $\eta_{\xi}$ is
precisely $BZ$. Hence
\ben
\eta_{\xi} & \in & Fol(\dim D,\dim D-(n-k), \dim BZ,0) \\
\eta_{\xi} & \in & Fol(2n-k-r-1, n-r-1, n-k-1, 0).
\een
Note that the set of foliation germs to which $\eta_{\xi}$ belongs is
of the form $Fol(n-s,k-s,r-s,0)$.
\par
Next $\eta_{\xi}$ can further be decomposed. Each slice
\[ \{ (\ast ,...,\ast ,\lambda_2,...,\lambda_{n-k},0,...,0)
                                                            \in\C^n\} \]
of $D$ is foliated by the single equation
$X_1^2+\cdots +X_{n-r}^2=\lambda_1$.
Let this foliation be $\alpha$. Then
\ben
\alpha & \in & Fol(\dim \mbox{\rm ~of~slice,~(dim~of~slice)}-1, 0,0) \\
\alpha & \in & Fol(n-r,n-r-1,0,0).
\een
Hence finally $\eta_{\xi}=i_{n-k-1}(\alpha )$ and $\xi =\eta_{\xi}\times
\C^{k+r+1-n}$, showing  that $\xi$ is a split foliation.
\atla
{\bf Remarks } \par
The problem discussed here is closely related to the classical extension
problem of coherent sheaves, see \cite{SE} and \cite{SI}.
For further geometry of singular
holomorphic foliations, and their residues,
see \cite{SS} and \cite{SU1}. \\
Sinan Sert\"{o}z \\
Bilkent University \\
Mathematics Department \\
06533 Ankara, Turkey \atla
e-mail: sertoz@fen.bilkent.edu.tr

\begin{thebibliography}{99}
%
\bibitem{AV} J.C. Alexander \& A. Verjovsky,
First integrals for singular
holomorphic foliations with leaves of bounded volume, Springer Lecture
Notes in Mathematics no: 1345, (1986) 1-10.
%
\bibitem{BA} P. Baum, Structure of foliation singularities, Advances in
Math., 15 (1975) 361-374.
%
\bibitem{BB} P. Baum \& R. Bott,
Singularities of holomorphic foliations,
J. of Differential Geometry, 7 (1972) 279-342.
%
\bibitem{CA} C. Camacho \& A.L. Neto, Geometric Theory of Foliations,
Birk\"{a}user, 1985.
%
\bibitem{CE} B. Cenkl,
Residues of singularities of holomorphic foliations,
J. of Differential Geometry, 13 (1978) 11-23.
%
\bibitem{GR} P.A. Griffiths, Two theorems on extensions of holomorphic
mappings, Inventiones Math., 14 (1971) 27-62.
%
\bibitem{XA} X. Gomez-Mont,
The transverse dynamics of a holomorphic flow,
Annals of Math., 126 (1987).
%
\bibitem{GSV} X. Gomez-Mont, J. Seade \& A. Verjovsky (Eds.),
Holomorphic Dynamics, Proceedings, Mexico 1986, Springer Lecture
Notes in Mathematics no: 1345.
%
\bibitem{MA} B. Malgrange, Frobenius avec singularites, Inventiones
Math., 39 (1977) 67-89.
%
\bibitem{SE2} J-P. Serre, Faisceaux alg\'{e}briques coh\'{e}rents,
Annals of Math., 61 (1955) 197-278.
%
\bibitem{SE} J-P. Serre, Prolongement de faisceaux analytique
coh\'{e}rent, Annales de L'Institu Fourier, Grenoble, 16 (1966) 363-374.
%
\bibitem{SS} S.Sert\"{o}z, Residues of singular holomorphic foliations,
Compositio Math., 70 (1989) 227-243.
%
\bibitem{SI} Y. Siu \& G. Trautmann,
Gap sheaves and extension of coherent analytic subsheaves,
Springer Lecture Notes in Mathematics, no: 172, (1971)
%
\bibitem{SU1} T. Suwa, Singularities of complex analytic foliations,
Proc. of Sympos. in Pure Math., vol 2 part 2 (1983) 551.
%
\bibitem{SU2} T. Suwa, Determinacy of analytic foliation germs, Studies
in Pure Math. 5, Foliations, Kinokuniye Co. Ltd., (1985) 427-460.
%
\bibitem{TH1} R. Thom, On singularities of foliations, Manifolds-Tokyo
(1973) 171-174.
%
\bibitem{TH2} R. Thom, Sur les bouts d'une feuilee d'un feuilletage au
voisinage d'un point singulier isole, Springer Lecture Notes in
Mathematics no: 1345, (1986) 317-321.
%
\bibitem{VI} S.M. Visik, Singularities of analytic foliations and
characteristic classes, Functional Anal. Appl. 7 (1973) 1-15.
%
\end{thebibliography}
\end{document}